\documentclass{gtart}
\usepackage{amssymb,rlepsf}
\usepackage[dvips]{graphicx}
\usepackage{psfrag}

\input gtmonout
\volumenumber{1}
\volumeyear{1998}
\volumename{The Epstein birthday schrift}
\pagenumbers{117}{125}
\received{13 May 1998}
\revised{15 October 1998}
\published{21 October 1998}
\papernumber{5}

\newtheorem{lem}{Lemma}
\newtheorem{tem}{Theorem}
\def\<{\langle}
\def\>{\rangle}

\def\adjpsfrag <#1,#2> #3#4{\psfrag{#3}%
{\smash{\rlap{\kern #1 \raise #2\hbox{#4}}}}}

\begin{document}

\title{All Fuchsian Schottky groups are\\classical Schottky groups}
\shorttitle{Fuchsian Schottky groups}
\author{Jack Button}

\address{Wadham College, Oxford, OX1 3PN, UK}
\email{button@maths.ox.ac.uk}

\begin{abstract}
Not all Schottky groups of M\"{o}bius transformations are classical
Schottky groups. In this paper we show that all Fuchsian Schottky
groups are classical Schottky groups, but not necessarily on the
same set of generators.
\end{abstract}
\asciiabstract{%
Not all Schottky groups of Moebius transformations are classical
Schottky groups. In this paper we show that all Fuchsian Schottky
groups are classical Schottky groups, but not necessarily on the
same set of generators.}

\maketitle

\section{Introduction}

A Schottky group of genus $g$ is a group of M\"{o}bius transformations
acting on the Riemann sphere $\overline{\mathbb C}$ generated by $g$ elements
$A_i,1\leq i \leq g$, each of which possesses a pair of Jordan curves
$C_i,C'_i\subseteq \overline{\mathbb C}$, with the property that the $2g$
curves are mutually disjoint and that $A_i$ maps $C_i$ onto $C'_i$ where
the outside of $C_i$ is sent onto the inside of $C'_i$. Direct use of
combination theorems tells us that the resulting group is free on $g$
generators, is discrete with a fundamental domain the region exterior
to the $2g$ curves, and consists entirely of loxodromic and hyperbolic
elements.

If in addition we can take all the Jordan curves to be geometric circles
then the resulting group is called a classical Schottky group (or sometimes
in order to be more specific we say it is classical on the generators
$A_1,\ldots,A_g$). Marden \cite{Mar74} showed that not all Schottky groups
are classical Schottky groups. Put very briefly, he argued that the algebraic
limit of classical Schottky groups must be geometrically finite and so his
isomorphism theorem implies that the ordinary set $\Omega$ of this limit
cannot be empty. But most groups on the boundary of Schottky space have
an empty ordinary set, so Schottky space strictly contains classical
Schottky space. However, this argument is certainly non-constructive,
raising the question of finding an explicit nonclassical Schottky group.
Zarrow \cite{Z75} claimed to have found such an example, but the paper of
Sato \cite{S88} shows that it is in fact a classical Schottky group.
A little later Yamamoto \cite{Y91} did construct a nonclassical
Schottky group.

The purpose of this paper is to show that if we examine the most 
straightforward cases where we might expect to find a counterexample,
namely Fuchsian Schottky groups, then this approach is doomed to failure
as all such groups are classical Schottky groups.
Specifically we show that:

({\bf 1})\qua Given a Fuchsian Schottky group $G$ of any genus $g$ then there
exists a generating set for $G$ of $g$ hyperbolic M\"{o}bius
transformations on which $G$ is classical.\\

({\bf 2})\qua  The Fuchsian Schottky group $G$ is classical on all possible
generating sets if and only if $g=2$ and $G$ is generated by
a pair of hyperbolic elements with intersecting axes.\\

({\bf 3})\qua  There exists a  Fuchsian group which is Schottky on a particular
generating set, but which cannot be classical on those generators.\\
\ppar

The author would like to thank the referee for comments on an earlier
draft of this paper.

\section{Proof of Main Theorem}

Given any finitely generated Fuchsian group $G$ (namely a discrete
subgroup of $PSL(2,\mathbb R)$) containing no elliptic elements, we
form the quotient surface $S=U/G$ where $U$ is the upper half
plane. The complete hyperbolic surface $S$ has ideal boundary
$\partial S=(\overline{\mathbb R}\cap\Omega_G)/G$, where
$\overline{\mathbb R}$ is the boundary of $U$ in the Riemann sphere
$\overline{\mathbb C}$ and $\Omega_G$ is the ordinary set of $G$.
Note that $G$ is Schottky if and only if $S$ is a closed surface minus
at least one hole (although $S$ cannot be a one-holed sphere).  This
is because a Fuchsian group $G$ with a quotient surface $S$ as above
must be free and purely hyperbolic, and this implies (see, say
\cite{Mas?Mar24}) that $G$ is indeed Schottky.

If $S$ is a surface of genus $n$ with $h$ holes then $G$ will be a free
group of some rank $r$. The process of doubling $S$ along its boundary
corresponds to considering the quotient of the whole ordinary set
$\Omega_G$ by $G$. As $G$ is a Schottky group, $\Omega_G/G$ is topologically
a closed surface of genus $r$. Therefore we conclude that $r=2n+h-1$
(with $n\geq 0,h\geq 1$ and $r\geq 1$).

The idea of the proof of theorem 1
is that given any such surface $S=U/G$, we find a particular reference
surface, homeomorphic to $S$, which has
a system of simple closed geodesics $\gamma_1 ,\ldots ,\gamma_r$ 
corresponding  to a generating set for $G$. We also find  disjoint
complete simple geodesics $l_1 ,\ldots , l_r$ 
on this reference surface which are properly
embedded (they can be thought of as having their endpoints up
the ``spouts''), where
$l_i$ intersects $\gamma_i$ once and is disjoint from $\gamma_j\ (j \neq i)$.
We will find that
if we cut  along these geodesics $l_1, \ldots ,l_r$, a disc is obtained.
We are then able to transfer these curves across to $S$.
By viewing the process upstairs in the upper half plane $U$ we get
a fundamental domain for $G$, and then
we can see directly that $G$ is classical Schottky on our generating set.

\begin{tem}
Given a Fuchsian Schottky group $G$ of any genus $g$ then there
exists a generating set for $G$ of $g$ hyperbolic M\"{o}bius
transformations on which $G$ is classical.
\end{tem}
\begin{proof} We prove the result by taking a standard Fuchsian classical
Schottky group $G_{n,h}$ for each possible topological surface of genus $n$
and $h$ holes, and transfer the two sets of geodesics to curves on any
other surface homeomorphic to $U/G_{n,h}$. These can be replaced by geodesics
with all necessary properties preserved.

First consider $h=1$. We choose $2n$ hyperbolic elements $A_1,\ldots ,A_{2n}$
so that their axes  all intersect at the same point, and ensure that
$G_{n,1}=\<A_1,\ldots ,A_{2n}\>$ is classical Schottky by choosing the
multipliers of the $A_i$ in order to obtain for each group $\<A_i\>$ a
fundamental domain $\Delta_i$ 
consisting of the intersection of the exteriors of
two geodesics $L_i$ and $L'_i=A_i(L_i)$ so that all conditions of the
free product combination theorem are satisfied; namely that
\[ \Delta_i\cup \Delta_j=U\mbox{ for }i\neq j\mbox{ and }
\bigcap_i\Delta_i\neq\emptyset.\] 
Then we have a fundamental
domain $\Delta_{n,1}$ (homeomorphic to a disc) for the discrete group
$G_{n,1}$.
There is one cycle of boundary intervals and so
by the discussion above, the surface $S_{n,1}=U/G_{n,1}$ is indeed of
genus $n$ with boundary a circle. 

We can project the axes of $A_i$ down onto the surface to obtain our simple
closed geodesics $\gamma_i$, and do the same with each $L_i$, which gives
us the complete simple geodesic $l_i$ right up to its two endpoints on
the boundary. These have the appropriate properties mentioned earlier, and
we see that the surface becomes a disc after cutting along all the
geodesics $l_1,\ldots l_{2n}$.

The group $G_{2,1}$ and the projection of these geodesics are illustrated
in figures 1 and 2.

\begin{figure}[htb]
\cl{\raise 2cm \hbox{\includegraphics[angle=-90,width=5cm]{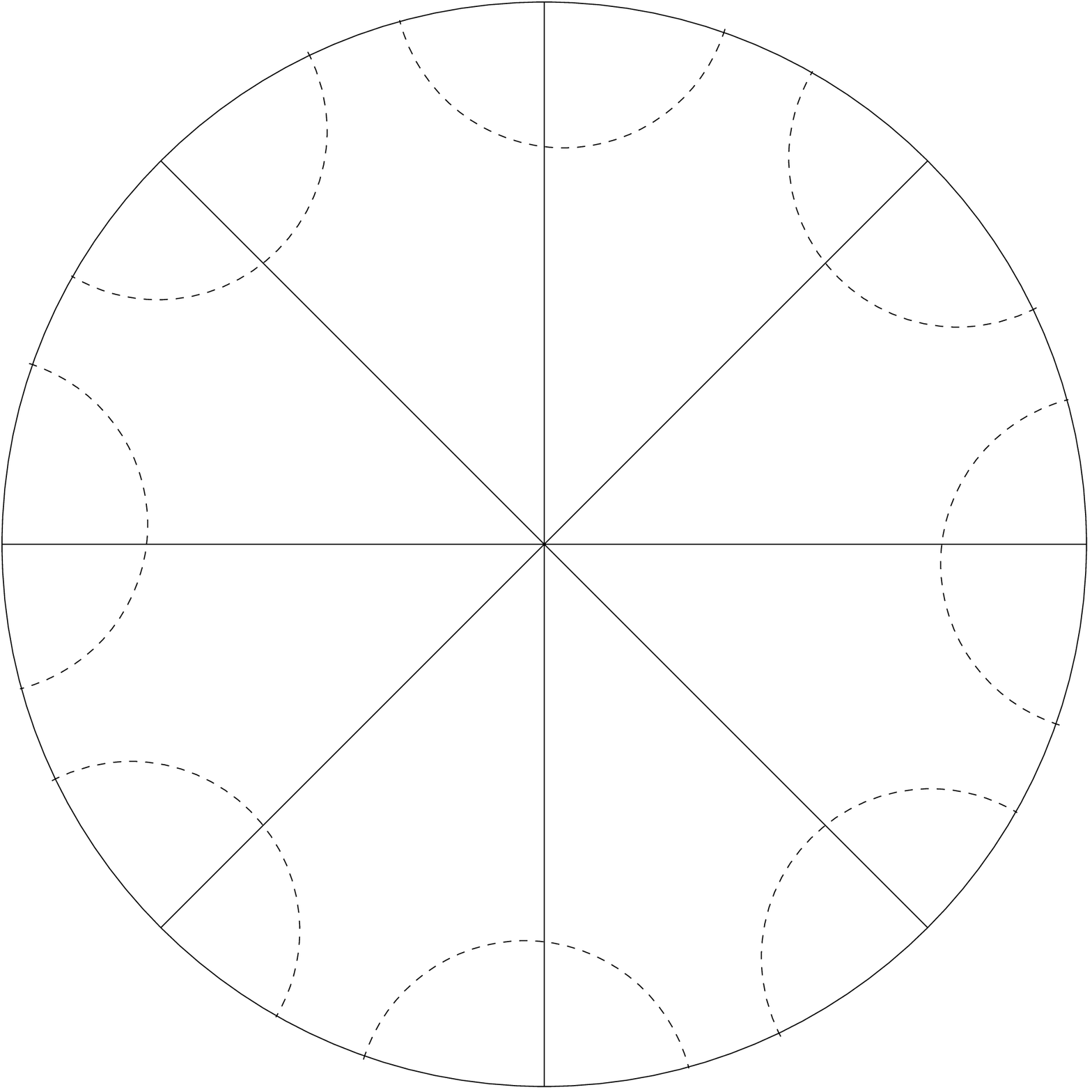}}
\hglue 1cm\lower 3cm \hbox{\includegraphics[width=3.5cm]{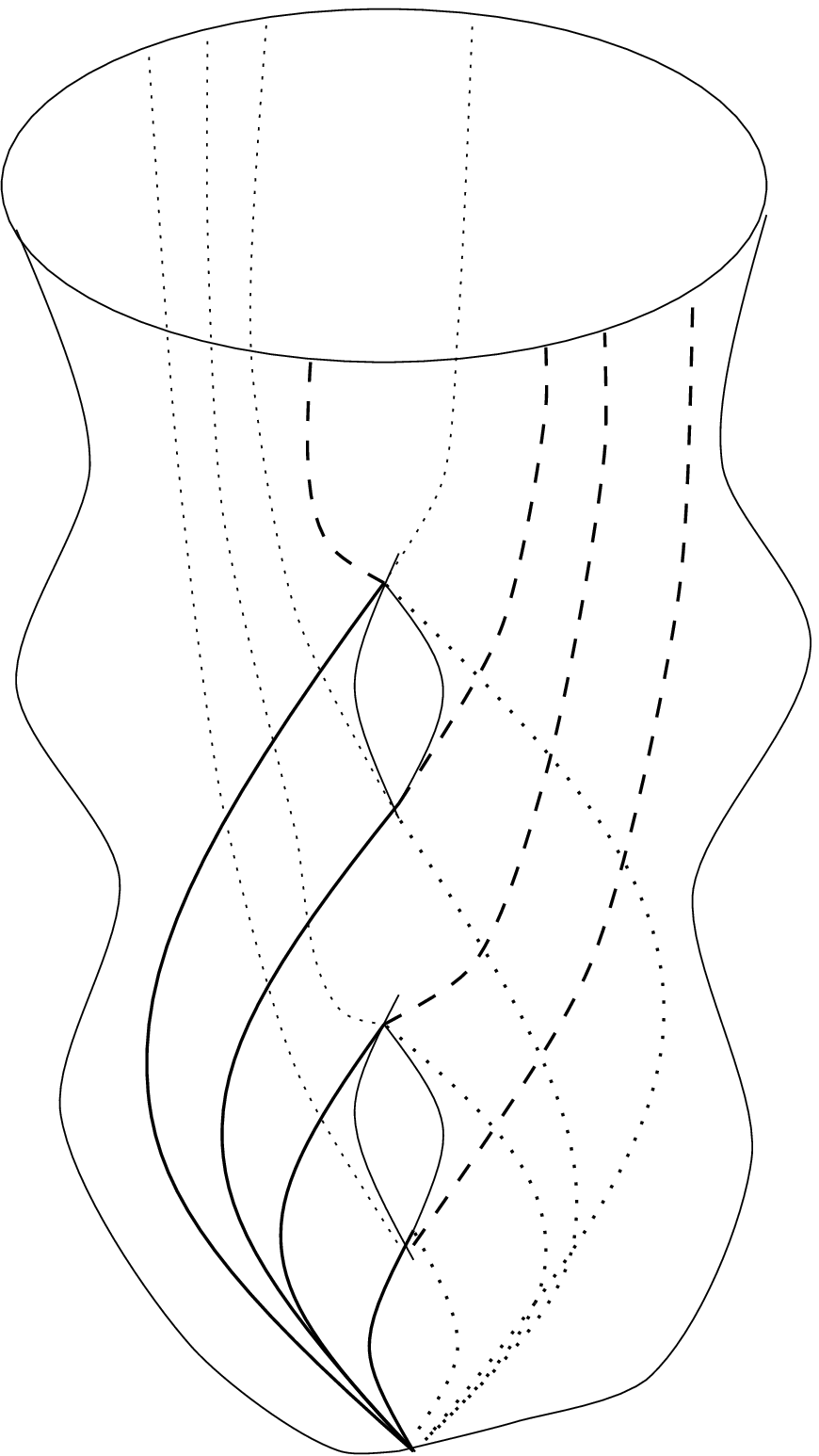}}}
\vglue 3mm
{\noindent\hbox{}\hglue 3.6cm Figure 1 \hglue 3.9cm Figure 2 \hfill\break}
\end{figure} 

In order to construct $G_{n,h}$ when $h\geq 2$, take $G_{n,1}$ and choose
an open interval $I$ between one endpoint of some $L_i$ and the nearest
endpoint of a neighbouring geodesic $L_j$. This interval lies inside the
ordinary set of $G_{n,1}$. Then inductively nest $h-1$ geodesics inside
the previous one, so that each geodesic has endpoints in $I$. We then find
hyperbolic transformations $A_{2n+1},\ldots ,A_{2n+h-1}$ with axes these
geodesics and with each transformation having two geodesics $L_i$ and
$L'_i=A_i(L_i)$, where $2n+1\leq i \leq 2n+h-1$, which it pairs. If these
fundamental domains are correctly placed then $G_{n,h}=\<A_1,\ldots A_{2n+h-1}
\>$ is a discrete group having the correct quotient surface
$S_{n,h}=U/G_{n,h}$ with a disc for a fundamental domain $\Delta_{n,h}$, where
$\partial \Delta_{n,h}$ consists of $4n+2h-2$ geodesics $L_i$ and $L'_i$,
along with the same number of intervals of $\overline{\mathbb R}$. The
geodesics and intervals alternate as we go round the boundary of the disc.
Also the projections of these axes and of these paired geodesics which
define $\gamma_i$ and $l_i$ have all the same properties as mentioned before.
The case $n=1$, $h=5$ is pictured in figures $3$ and $4$.

\begin{figure}[htb]
\cl{%
\includegraphics[angle=-90,width=6cm]{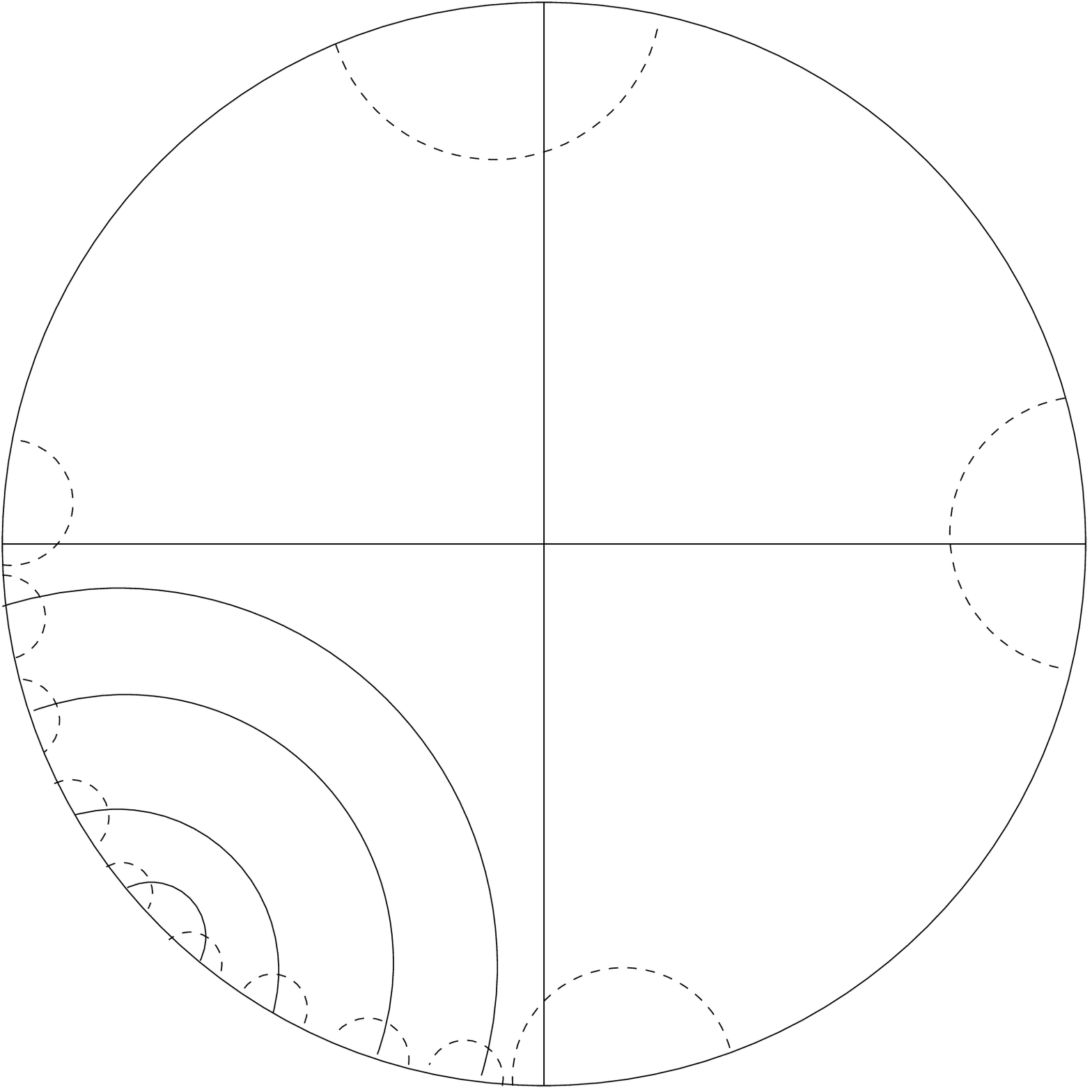}
\includegraphics[angle=-90,width=6cm]{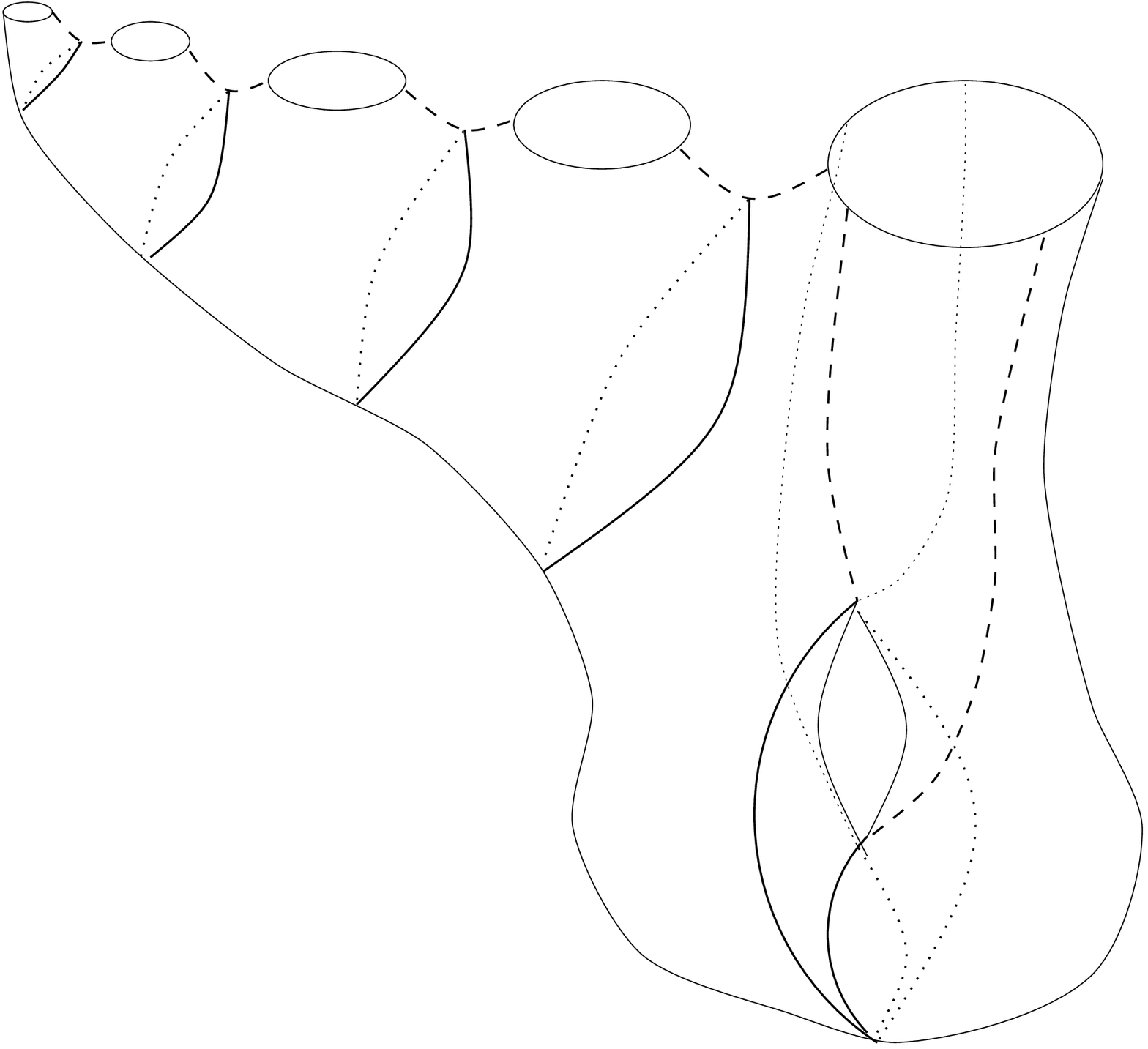}}
\vglue 3mm
{\noindent\hbox{}\hglue 2.8cm Figure 3 \hglue 6.2cm Figure 4 \hfill\break}
\end{figure}

Now given any Fuchsian Schottky group $G$ with quotient surface $S$ and
boundary $\partial S$, there exists a homeomorphism
\[ h\co S_{n,h}\cup \partial S_{n,h} \mapsto S \cup\partial S \]
for some $n$ and $h$. We also have natural continuous  projections
\[ \begin{array}{rlcl}
p\,\co  & U\cup (\Omega_{G_{n,h}} \cap \overline{\mathbb R}) & \mapsto &
S_{n,h}\cup \partial S_{n,h}\\
q\,\co  & U\cup (\Omega_{G} \cap \overline{\mathbb R}) & \mapsto & S\cup\partial S
\end{array} \]
where $p$ and $q$ are both covering maps, and both domains are simply
connected covering spaces of their images (where the elementary neighbourhoods
of points downstairs are open discs, or half discs for points on the
boundary).

By the lifting theorem, we have a continuous map
\[ H\co U \cup (\Omega_{G_{n,h}}\cap\overline{\mathbb R})\mapsto U\cup (\Omega_{G}
\cap \overline{\mathbb R}) \]
which is a lift of $hp$, so that $hp=qH$. By reversing $p$ and $q$, we see
that $H$ is a homeomorphism.

Take any element $g\in G_{n,h}$. This is a deck transformation of $p$ and so
$pg=p$. Conjugating $g$ by $H$, we have $q(HgH^{-1})=q$, thus $HgH^{-1}$ is
a deck transformation of $q$ and therefore $H$ defines an isomorphism of
$G_{n,h}$ onto $G$ by conjugation.

Note that $H$ maps $U$ to $U$ and $\Omega_{G_{n,h}}\cap \overline{\mathbb R}$
to $\Omega_G\cap \overline{\mathbb R}$, because it is a lift of $h$ which
sends boundary points to and from boundary points. Therefore the image
under $H$ of
the fundamental domain
$\Delta_{n,h}$ is a disc in $U$. But $H(\partial\Delta_{n,h})$ will consist of
$4n+2h-2$ disjoint closed intervals of $\overline{\mathbb R}$, along
with curves $H(L_i)$ and $H(L'_i)$ lying entirely in $U$ apart from their
endpoints
which
are also endpoints of these intervals of $\overline{\mathbb R}$. We find
that the order in which the images under H of the $L_i$, $L'_i$ and 
the intervals
appear around $\partial H(\Delta_{n,h})=H(\partial\Delta_{n,h})
\subseteq U\cup (\Omega_G
\cap\overline{\mathbb R})$ is the same as the original order around
$\partial\Delta_{n,h}$ (or the opposite order if $H$ is orientation
reversing).

By setting  $B_i=HA_iH^{-1}$ we obtain a generating set for $G$, and
because $A_i$ sends the geodesic $L_i$ to $L'_i$, we see that $B_i$
sends the curve $H(L_i)$ to the curve $H(L'_i)$. Also it is easy to
check that the disc $H(\Delta_{n,h})$ is a fundamental domain for the
action of $G$ on $U$. In particular, the intersection of the exteriors
in $U$ of $H(L_i)$ and $H(L'_i)$ is a fundamental domain for $\<B_i\>$.
We replace these two curves by geodesics $M_i$ and $M'_i=B_i(M_i)$ which
have the same endpoints. Just as in \cite{DN}, this gives us $2n+h-1$
pairs of geodesics freely homotopic to the curves they replaced, and paired
by a generating set ${B_i}$ with another fundamental domain $D_i$
for each group   
$\<B_i\>$ that lies between these two geodesics. The free product
combination theorem can be applied to $\<B_1\>,\ldots \<B_{2n+h-1}\>$, as 
$D_i\cup D_j=U$ for $i\neq j$ and $\bigcap_i D_i\neq\emptyset$.
We can see this by looking at the endpoints of the geodesics which have
not been changed when passing from curves. Therefore, by reflecting this
picture in the real axis, the group $G$ is generated by elements $B_i$,
each of which possesses a pair of mutually disjoint geometric circles
$C_i$ and $C'_i$, with the outside of $C_i$ being sent by $B_i$ onto
the inside of $C'_i$. By definition, $G$ is a classical Schottky group.
\end{proof}

\begin{figure}[htb]\small
\adjpsfrag <-5pt,-4pt> {A}{$A$}
\adjpsfrag <0pt,-4pt> {B}{$B$}
\adjpsfrag <1pt,-3pt> {AB}{$AB$}
\adjpsfrag <-3pt,-3pt> {BA}{$BA$}
\adjpsfrag <-10pt,13pt> {B1A}{$B^{-1}A$}
\adjpsfrag <-6pt,-1pt> {BA1}{$BA^{-1}$}
\begin{center}
\includegraphics[width=6.5cm,angle=-90]{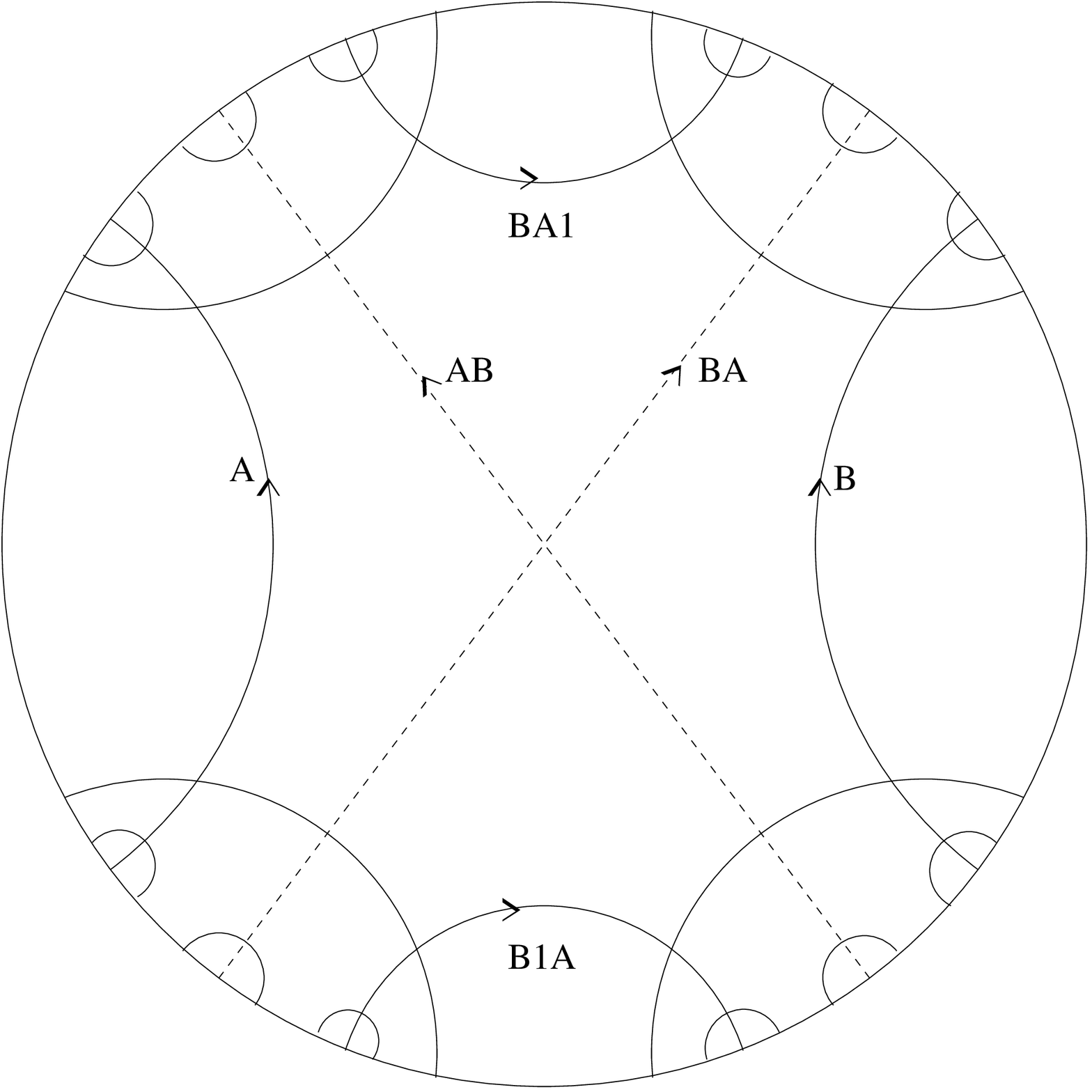}
\vglue 3mm
Figure 5
\end{center}
\end{figure}

\section{Proof of other Theorems}

Suppose we are given any two hyperbolic elements $A$ and $B$ with different
axes. We want to know when $G=\<A,B\>$ is free, discrete and purely
hyperbolic (hence Schottky). This problem falls naturally into two cases.\\

\begin{figure}[htb]\small
\adjpsfrag <10pt,-4pt> {A}{$A$}
\adjpsfrag <-14pt,-4pt> {B}{$B$}
\adjpsfrag <-2pt,-2pt> {Ax}{$Ax$}
\adjpsfrag <-6pt,-3pt> {Bx}{$Bx$}
\adjpsfrag <-6pt,15pt> {B1A}{$B^{-1}A$}
\adjpsfrag <-4pt,-1pt> {BA1}{$BA^{-1}$}
\adjpsfrag <-1pt,-12pt> {x}{$x$}
\adjpsfrag <-1pt,-0pt> {y}{$y$}
\adjpsfrag <-2pt,-0pt> {z}{$z$}
\adjpsfrag <-13pt,-15pt> {Ay}{$A^{-1}y$}
\adjpsfrag <-0pt,-15pt> {Bz}{$B^{-1}z$}
\begin{center}
\includegraphics[width=7cm,angle=-90]{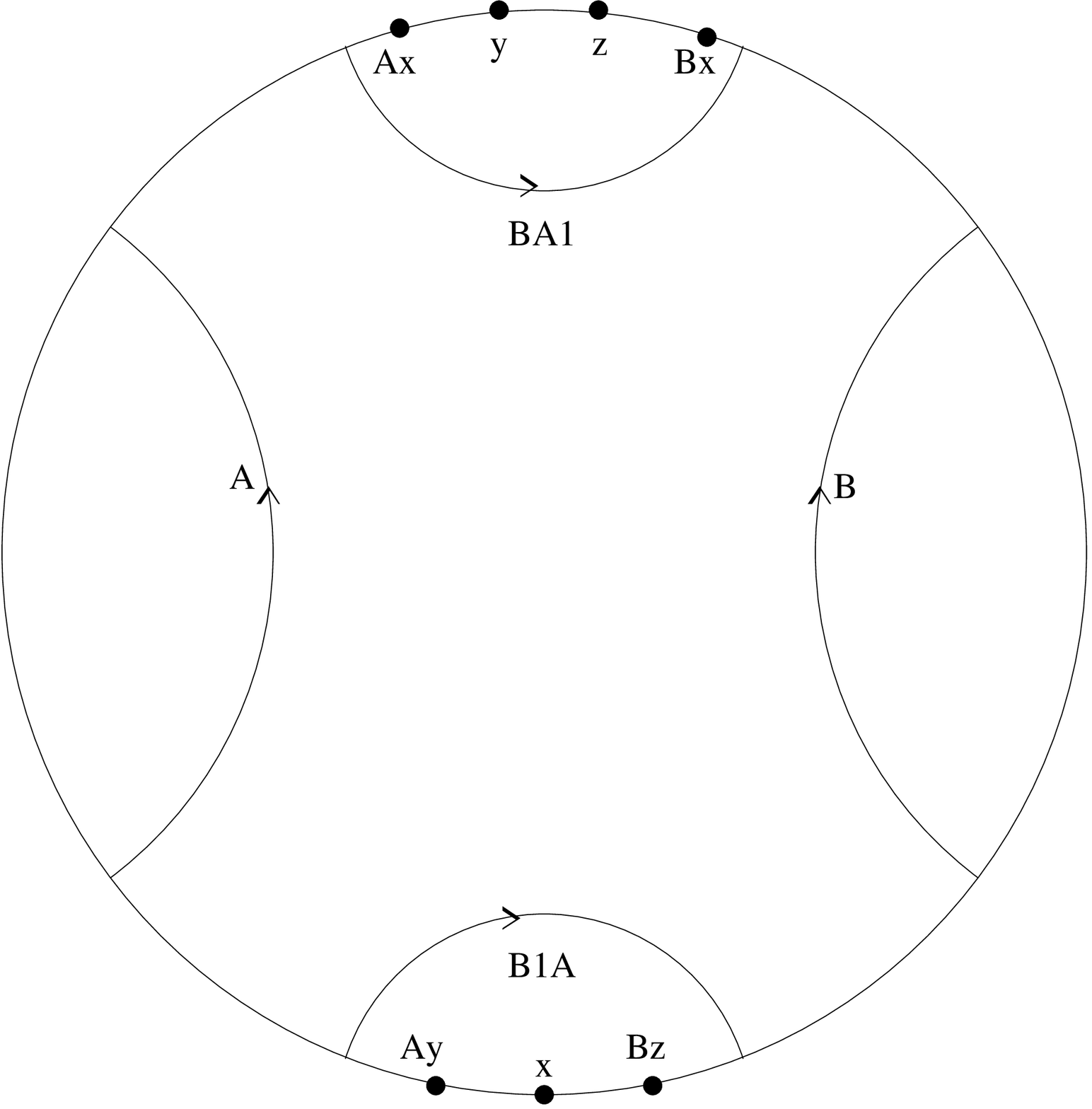}
\vglue 6mm
Figure 6
\end{center}
\end{figure}
({\bf A})\qua  The two hyperbolic elements have intersecting axes. Then it is
well known that $G$ is free, discrete and purely hyperbolic if and only if
the commutator $ABA^{-1}B^{-1}$ is hyperbolic. See for instance \cite{P72}
where this is shown by explicitly exhibiting two pairs of geometric circles,
one paired by $A$ and one by $B$. In this case the quotient surface is a one
holed torus and, as any generating pair will have intersecting axes, we see
that $G$ is classical on every possible generating pair.

Alternatively we can see this directly from section 1 by using the fact that
there will exist a homeomorphism from our standard surface to the quotient
surface of $G$ that takes the two simple closed geodesics $\gamma_1,\gamma_2$
onto two curves freely homotopic to the
simple closed geodesics corresponding to any generating
pair of $G$.\\

({\bf B})\qua  The hyperbolic elements have non-intersecting axes. If so
then all generating pairs of $G$ must have non-intersecting axes, or
else we are back in case (A).

First suppose $G$ is a classical Schottky group on these two generators $A$
and $B$. Without loss of generality we can replace any generator by its
inverse so that we get a picture such as the one in figure 5, with the
arrows on the two generators in the same direction. The quotient surface
is a three holed sphere. Note that the axis of $AB$ projects down onto a
``figure of eight'' geodesic, and so this group cannot be classical on the 
generating pair $\<A,AB\>$. 

\begin{tem}
A group $G$ that has a quotient surface which is not a one holed torus cannot
be classical on all generating sets.
\end{tem}

\begin{proof}
We have already considered any $G$ generated by two elements.
Given any $G$ generated by three or more elements, we can find a pair of
generators with non-intersecting axes, and use the above argument on the
subgroup generated by this pair. As the subgroup is not classical on all
generating sets, nor is $G$.
\end{proof}

Finally we show the existence of a Fuchsian group generated by two elements
which is Schottky, but not classical, on this generating pair.

\begin{lem}
A group $G=\<A,B\>$ (where $A$ and $B$ are hyperbolic elements with
non-intersecting axes, 
oriented as in figure 5)
is classical
on $\<A,B\>$ if and only if both fixed points of $B^{-1}A$  lie in the
interval between the repelling fixed points of $A$ and $B$.
\end{lem}

\begin{proof}
If we know $G$ is classical on $\<A,B\>$ then we can build up a pattern of
nested circles as in figure 5, and see the location of the fixed points
of the axes directly. Conversely if we only have information as in figure
6 then we consider the image of a suitable point $x$ under the generators.

The axis of $B^{-1}A$ is sent to the axis of $BA^{-1}$ by both generators,
and also note that the arrows on $BA^{-1}$ and $B^{-1}A$ are as in the picture
(for instance consider the image of a fixed point of $A$).
Then we choose any $x$ inside the interval enclosed by the axis of $B^{-1}A$,
and mark it and its images under $A$ and $B$. We can take any two points
$y$ and $z$ in the interval between $Ax$ and $Bx$, and use these as endpoints
for the geometric circles we require.

We can see that $A^{-1}y$ will be closer than $x$
to the repelling fixed point of $A$, and similarly with $B^{-1}z$ and $B$.
This gives us four endpoints $y,z,A^{-1}y$ and $B^{-1}z$, one for each circle.
We have four more
endpoints to mark but this choice is totally arbitrary: merely pick any
point in the interval between $A$'s fixed points, along with its image
under $A$, and do the same for $B$ too. This provides us with our two
pairs of circles which show that $G$ is discrete, and classical on $\<A,B\>$.
\end{proof}
\begin{figure}[htb]\small
\adjpsfrag <-5pt,10pt> {A}{$A$}
\adjpsfrag <-0pt,-2pt> {B}{$B$}
\adjpsfrag <0pt,0pt> {C}{$C$}
\adjpsfrag <-2pt,-1pt> {CA}{$C_A$}
\adjpsfrag <0pt,-4pt> {CB}{$C_B$}
\adjpsfrag <-1pt,2pt> {CA1}{$C'_A$}
\adjpsfrag <1pt,1pt> {CB1}{$C'_B$}
\begin{center}
\includegraphics[width=6cm,angle=-90]{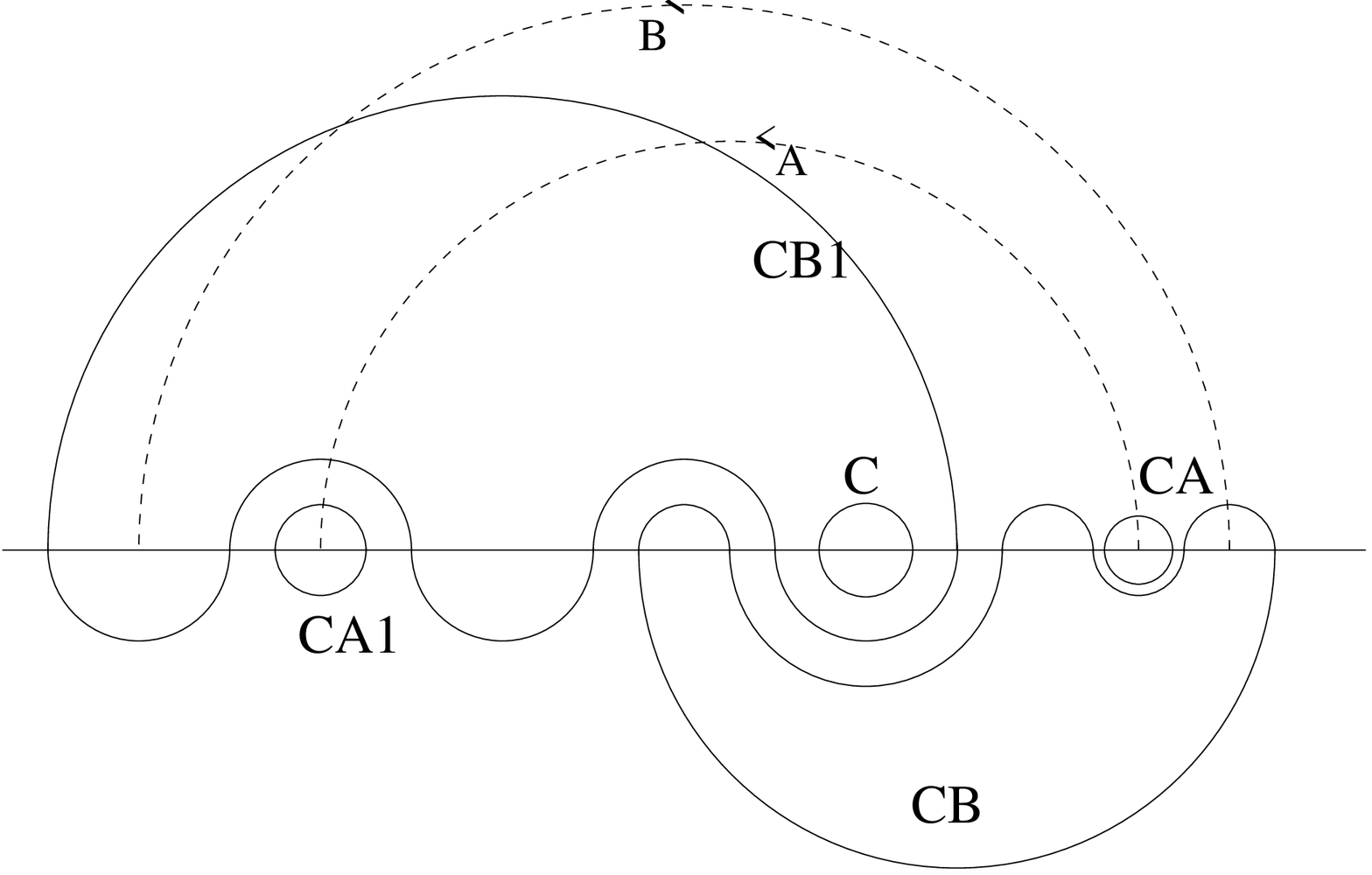}
\vglue 2mm
Figure 7
\end{center}
\end{figure}
\begin{tem}
The Fuchsian group in figure 7, which is Schottky on the generators $A$ and
$B$, is not classical on them.
\end{tem}
\begin{proof}
The exterior $F$ of the two pairs of curves $C_A,C'_A$ (paired by $A$)
and $C_B,C'_B$ (paired by $B$)
is a fundamental domain, and is sent by the
element $BA^{-1}$ inside the circle $C(=B\left(C_A\right) )$. The attracting
fixed point of $BA^{-1}$ must lie inside $C$ and therefore it separates
the fixed points of $A$.
\end{proof}

\Addresses

\end{document}